\documentclass[english]{IEEEtran}
\usepackage[T1]{fontenc}
\usepackage[latin9]{inputenc}
\usepackage{color}
\usepackage{amssymb}
\usepackage{setspace}
\usepackage{color}
\usepackage{mathtools}
\usepackage{algorithmicx}
\usepackage[titlenumbered,ruled]{algorithm2e}
\DeclarePairedDelimiter{\ceil}{\lceil}{\rceil}
\makeatletter
\usepackage{balance}
\usepackage{babel}
\usepackage{multirow}
\usepackage{amsmath}
\usepackage{graphicx}
\usepackage{bigstrut}
\renewcommand{\fnum@figure}{Fig.~\thefigure} 
\usepackage{cite}
\usepackage{graphics}
\usepackage{adjustbox}
\usepackage{color}
\usepackage{soul}
\usepackage{url}
\newtheorem{example}{Example}[section]
\allowdisplaybreaks
\makeatother

\usepackage{babel}
\begin{document}

\title{AC Transmission Network Expansion Planning: A Semidefinite Programming
Branch-and-Cut Approach}

\author{Bissan Ghaddar and Rabih A. Jabr, \emph{Fellow}, \emph{IEEE} \thanks{B. Ghaddar is with the Department of Management Sciences, University
of Waterloo, Canada (email: bissan.ghaddar@uwaterloo.ca). }\thanks{R. A. Jabr is with the Department of Electrical \& Computer Engineering,
American University of Beirut, P.O. Box 11-0236, Riad El- Solh / Beirut
1107 2020, Lebanon (email: rabih.jabr@aub.edu.lb).}}
\maketitle
\begin{abstract}
Transmission network expansion planning is a mixed-integer optimization
problem, whose solution is used to guide future investment in transmission
equipment. An approach is presented to find the global solution of
the transmission planning problem using an AC network model. The approach
builds on the semidefinite relaxation of the AC optimal power flow
problem (ACOPF); its computational engine is a new specialized branch-and-cut
algorithm for transmission expansion planning to deal with the underlying mixed-integer ACOPF problem. Valid inequalities
that are based on specific knowledge of the expansion problem are
employed to improve the solution quality at any node of the search
tree, and thus significantly reduce the overall computational effort
of the branch-and-bound algorithm. Additionally, sparsity of the semidefinite relaxation is exploited to further reduce the computation time at each node of the branch-and-cut tree.
Despite the vast number of publications on transmission expansion
planning, the proposed approach is the first to provide expansion
plans that are globally optimal using a solution approach for the mixed-integer ACOPF problem.
The results on standard networks serve as important benchmarks to
assess the solution quality from existing techniques and simplified
models.
\end{abstract}

\begin{IEEEkeywords}
Design optimization, mathematical programming, nonlinear network analysis,
optimization methods, power system planning.
\end{IEEEkeywords}

\section{Introduction}

\IEEEPARstart{T}{he} optimal transmission network expansion planning
(TNEP) problem seeks to compute the specifications and locations of
electrical transmission equipment so that the expanded system can
meet the expected future load and generation patterns \cite{Garver_1970}.
The planning problem can be carried out over a single stage, or it
can be extended over multiple stages that represent a longer planning
period; the optimal decision in multistage planning links the required
investment in transmission equipment to a particular time period \cite{Escobar_2004,Vinasco_2011}.
In its most accurate representation, TNEP is a mixed-integer nonlinear
program whose solution is computationally challenging. The application
of computing techniques to solve the TNEP problem has been steadily
increasing since the 1970s, and has passed through various accuracy
levels of network representation that were adapted to the existing
optimization solver technology \cite{Romero_2002}. This paper considers
the solution of TNEP with the system being accurately represented
via the AC network model; it employs the power network semidefinite-based
relaxation \cite{Lavaei_2012} in a branch-and-cut optimization methodology
that is specifically tailored to the TNEP problem. The outcome is
a TNEP plan employing the AC network model, that is provably the globally
optimal solution. It is worth noting that none of the existing AC-TNEP
techniques can claim the global optimally of their solutions.

The TNEP literature reports expansion solutions using different network
models: (i) the transportation model that enforces Kirchhoff's current
law (KCL) at all nodes \cite{Garver_1970,Haffner_2000,Haffner_2001},
(ii) the hybrid model that adds to the transportation model the Kirchhoff's
voltage law (KVL) equations only for existing circuits \cite{Villasana_1985},
(iii) the DC network model that adds to the hybrid model the KVL equations
for new circuits \cite{Romero_2007}, and (iv) the AC network model
\cite{Rider_2007}. Both the DC and AC network models when employed
in TNEP result in terms that have products of integer and continuous
variables, which further contribute to the challenge in computing
the solution. This can be circumvented by the use of the disjunctive
reformulation in TNEP problems that are based on DC \cite{Bahiense_2000}
and AC \cite{Zhang_2012_c} power flow models. 

The solution techniques for TNEP can be classified into three categories:
(i) the locally optimal search methods, (ii) the meta-heuristics,
and (iii) the complete search methods based on conventional optimization
theory. The constructive heuristic algorithm \cite{Garver_1970,Rider_2007}
is one of the earliest local search methods applied to TNEP; it can
give good quality solutions but which are unlikely to be optimal.
The discrepancy-bounded local search \cite{Bent_2012} is a generalization
of the constructive heuristic that gives better solutions with ACOPF models. The meta-heuristic methods are rooted in the processes
of natural and physical systems \cite{Gallego_1998a}; these include
simulated annealing, tabu search, and several variants of genetic
algorithms that mimic the mechanisms of natural selection and evolution
\cite{Gallego_1998b,da_Silva_2000,da_Silva_2001,Escobar_2004,de_J_Silva_2005,Rodriguez_2008,Rahmani20101056}.
Although meta-heuristics can avoid being trapped in locally optimal
TNEP solutions, they do not give a rigorous indicator on the quality
of the solution. The classical optimization algorithms tend to fall
into two categories: (i) Benders decomposition type algorithms and
(ii) Branch-and-Bound/Branch-and-Cut. Due to the non-convexity of
the DC/AC-TNEP formulations, it is possible for the Benders cut to
chop off a part of the feasible region that includes the global solution;
for DC-TNEP, this can be circumvented either via the hierarchical
decomposition method \cite{Romero_1994a,Romero_1994b} or by reformulating
the TNEP problem as a disjunctive model \cite{Binato_2001,AlizadehMousavi2016275}.
State-of-the-art implementations of branch-and-cut approaches are
currently provided in commercial packages for mixed-integer linear
programming; these have been used in disjunctive models of DC-TNEP
\cite{Vinasco_2011}, in addition to enhanced linear representations
that approximate network losses and reactive power constraints \cite{Zhang_2012,Zhang_2013,Bent_2014,Macedo_2016};
a recent extension includes modeling the choice of both HVAC and HVDC
transmission equipment \cite{Dominguez_2017}.

While the linear network models have been extensively used in TNEP,
the strict AC model is rarely discussed due to the computational complexity
that it introduces; only \cite{Rider_2007,Rodriguez_2008,Bent_2012,Rahmani20101056,Zhang_2012_c}
from amongst the above surveyed methods use the AC power flow equations.
However, TNEP solutions that are based on the DC network model can
yield arbitrarily poor designs that require further reinforcement
to satisfy the AC network feasibility \cite{Bent_2012,Bent_2014}.
The main impediment to the use of the AC network model in TNEP based
on classical optimization is the non-convexity of the AC power flow
constraints; this has been circumvented in \cite{Jabr_2013} by exploiting
(i) the conic relaxation of the AC optimal power flow (OPF) over networks
with virtual controllable phase shifters, and (ii) a state-of-the-art
branch-and-cut algorithm for mixed-integer conic programming. Although
\cite{Jabr_2013} reports solutions that satisfy the AC network model
constraints, these solutions may not be globally optimal due to the
use of linear equations that are needed to alleviate the effect of
the virtual phase shifters. To ascertain global optimality of the
AC-TNEP solutions, this paper proposes the use of the semidefinite
programming (SDP) relaxation of the AC optimal power flow problem
\cite{Lavaei_2012}; the SDP relaxation does not entail the disadvantage
associated with the virtual phase shifters, but it requires the development
of a specialized branch-and-cut solver to handle SDP constraints.
This paper contributes an SDP branch-and-cut method for AC-TNEP,
provides details of valid inequalities, and proposes to utilize sparsity and problem structure to speed up the solution. The
results are provably global, and serve as benchmarks against which
the quality of approximated methods and heuristic solutions could
be gauged. 

\section{Mathematical Formulation}

In this section, the same notation is used as in \cite{Jabr_2013}. The topology of the power system $P = (N, E)$ is represented as an undirected graph, where each vertex $n \in N$ is called a ``bus'' and each edge $e \in E$ is called a ``branch'' linking buses to one another. The parameter $|N|$ denotes the number of buses and $|E|$ denotes the number of branches. Let $G \subseteq N$ be the set of generators and $E \subseteq N \times N$ be the set of all branches and $N(l)$ are buses adjacent to bus $l$. Let $S^d_l=P^d_l + jQ^d_l$ be the active and reactive load (demand) at each bus $l \in N$ and $P^g_l + jQ^g_l$ represent the complex power of the generator at bus $l \in G$. Define $V_l =\Re{V_l}+ j \Im{V_l}$ as the voltage at each bus $l \in N$ and $S_{lm}=P_{lm}+jQ_{lm}$ as the apparent power flow on the line $(l,m) \in E$. The edge set $L \subseteq E$ contains the branches $ (l,m)$ such that the apparent power flow limit is less than a certain given tolerance $\varepsilon$. The next section describes the formulation of the transmission network expansion planning problem in details.

\subsection{SDP Formulation}

To formulate the TNEP model, the ACOPF problem is first presented as it is the main building block of the TNEP formulation. In this work, the focus is on the rectangular power-voltage formulation of the OPF problem. In the rectangular formulation, the bus voltages are represented by the real and imaginary voltage components. The real and reactive power
flows are quadratic functions of real and imaginary parts of the voltage. This results in the generator limits, the fixed loads and the apparent power line limits being non-convex constraints, and in addition the lower limits of bus voltage magnitudes are non-convex. Hence, the rectangular formulation is also non-convex. The following parameters are needed to formulate the problem:
\begin{itemize}
\item $P_l^{\min}$ and $P_l^{\max}$ are the limits on active generation capacity at bus $l$, where $P_l^{\min}=P_l^{\max} =0$ for all $l \in N / G$.
\item $Q_l^{\min}$ and $Q_l^{\max}$ are the limits on reactive generation capacity at bus $l$, where $Q_l^{\min}=Q_l^{\max} =0$ for all $l\in N / G$. 
\item $P_l^d$ is the active power demand at bus $l$.
\item $Q_l^d$ is the reactive power demand at bus $l$.
\item $V_l^{\min}$ and $V_l^{\max}$ are the limits on the voltage at a given bus $l$. 
\item $S_{lm}^{\max}$ is the limit on the absolute value of the apparent power of a branch $(l,m)\in L$.
\end{itemize}
Additionally, let $y_{lm}=g_{lm}+jb_{lm}$ be the series admittance in the $\pi$-model of line $(l,m)$ and $\bar{b}_{lm}$ be the shunt susceptance in the $\pi$-model of line $(l,m)$.

Given a complex voltage $V_l$ at bus $l$, let $\Re{V_l}$ denote the real part of $V_l$ and $\Im{V_l}$ denote the imaginary part. %The network model is constrained by the power flow equations:
% \begin{align}
% P^g_l &= P_l^d+\Re{V_l} \sum_{l=1}^n ( { \Re{y_{lm}} \Re{V_l} - \Im{y_{lm}} \Im{V_l} }) \notag \\
%     & + \Im{V_l} \sum_{l=1}^n ({ \Im{y_{lm}} \Re{V_l} - \Re{y_{lm}} \Im{V_l} }) \label{eqn:Pk} \\
% Q^g_l &= Q_l^d+\Re{V_l} \sum_{l=1}^n ({ - \Im{y_{lm}} \Re{V_l} - \Re{y_{lm}} \Im{V_l} }  ) \notag \\
%     & + \Im{V_l} \sum_{l=1}^n ({ \Re{y_{lm}} \Re{V_l} - \Im{y_{lm}} \Im{V_l} }) \label{eqn:Qk}
% \end{align}
In terms of power line flows, the following equations hold:
\begin{align}
P_{lm} &= b_{lm} ( \Re{V_l} \Im{V_m} - \Re{V_m} \Im{V_l}) \label{eqn:Plm}  \\
    &+ g_{lm}( \Re{V_l}^2 +\Im{V_m}^2- \Im{V_l} \Im{V_m}- \Re{V_l} \Re{V_m}) \notag \\
Q_{lm}  &= b_{lm} ( \Re{V_l} \Re{V_m} + \Im{V_l} \Im{V_m}-\Re{V_l}^2 -\Im{V_l}^2) \notag \\
    &+ g_{lm}( \Re{V_l}\Im{V_m}- \Re{V_m} \Im{V_l})  \notag \\
  & -\frac{\bar{b}_{lm}}{2}(\Re{V_l}^2 +\Im{V_l}^2) \label{eqn:Qlm}
\end{align}
Let $\Re{V}_l$ and $\Im{V}_l$ be the variables in addition to the variables $P^g_l$, $P_{lm}$ and $Q_{lm}$. The aim of the optimal power flow problem is to satisfy demand at all buses with the minimum total production costs of generators such that the
solution obeys the physical laws and other operational constraints such as transmission line flow limit constraints. The OPF can be formulated as the following quadratic non-convex problem: %\textcolor{red}{CHECK EQUATIONS}
\begin{align}
\tag{OPF} \label{eq:opf} \\
\min\ \ &\sum_{l\in G} c^2_l(P^g_l)^2+c^1_l(P^g_l)+c^0_l \label{eq1}\\
\text{s.t. }\ \ &P_l^{\min} \leq P_l^g \leq P_l^{\max}, \ \forall l \in G,  \label{eq2}\\
& Q_l^{\min}\leq Q^g_l\leq Q_l^{\max}, \ \forall l \in G, \label{eq3}\\
& P^g_l -P_l^d= \sum_{m\in N(l)} P_{lm},  \ \forall l \in N, \label{eq4} \\
& Q^g_l - Q_l^d= \sum_{m\in N(l)} Q_{lm},  \ \forall l \in N, \label{eq5}\\  
& (V_l^{\min})^2 \leq \Re{V}_l^2 +  \Im{V}_l^2 \leq  (V_l^{\max})^2,  \ \forall l \in N, \label{eq6}\\
& P_{lm}^2+Q_{lm}^2 \leq (S_{lm}^{\max})^2,  \ \forall (l,m)\in E,\label{eq7} \\
& \eqref{eqn:Plm}-\eqref{eqn:Qlm} \label{eq8}
% & P^g_l = P_l^d+\Re{V_l} \sum_{l=1}^n ( { \Re{y_{lm}} \Re{V_l} - \Im{y_{lm}} \Im{V_l} }) \notag \\
% & + \Im{V_l} \sum_{l=1}^n ({ \Im{y_{lm}} \Re{V_l} - \Re{y_{lm}} \Im{V_l} }) \label{eq6} \\
% & Q^g_l = Q_l^d+\Re{V_l} \sum_{l=1}^n ({ - \Im{y_{lm}} \Re{V_l} - \Re{y_{lm}} \Im{V_l} }  ) \notag \\
% & + \Im{V_l} \sum_{l=1}^n ({ \Re{y_{lm}} \Re{V_l} - \Im{y_{lm}} \Im{V_l} }) \label{eq7} \\
% &P_{lm} &= b_{lm} ( \Re{V_l} \Im{V_m} - \Re{V_m} \Im{V_l}) \label{eq8}  \\
%     &+ g_{lm}( \Re{V_l}^2 +\Im{V_m}^2- \Im{V_l} \Im{V_m}- \Re{V_l} \Re{V_m}) \notag \\
% & Q_{lm}  &= b_{lm} ( \Re{V_l} \Im{V_m} - \Im{V_l} \Im{V_m}-\Re{V_l}^2 -\Im{V_l}^2) \notag \\
%     &+ g_{lm}( \Re{V_l}\Im{V_m}- \Re{V_m} \Im{V_l}- \Re{V_m} \Im{V_l})  \notag \\
%   & -\frac{\bar{b}_{lm}}{2}(\Re{V_l}^2 +\Im{V_l}^2) \label{eq9}
\end{align}
The objective function \eqref{eq1} minimizes the cost of power generation. Constraints \eqref{eq2}, \eqref{eq3}, and \eqref{eq4} set limits on the active power, reactive power and represent the real and reactive nodal balance. Constraints \eqref{eq6} restricts the voltage on each bus. Constraints~\eqref{eq7} set a limit on the apparent power flow at every branch. While constraints~\eqref{eq8} define the power generated and the active and reactive power flow.    
To  apply  SDP  relaxation  to  the  rectangular  formulation  (OPF),  define  a  hermitian  matrix $X=VV^*$ that is $X=X^*$, where $X^*$ is the conjugate transpose of $X$. The Hermitian matrix $X$ is defined as:
$$X= \left[ \begin{matrix} |V_1|^2  & \cdots & V_1V_l^* & \cdots &  V_1V_n^*   \\
   \vdots  & \ddots & \vdots &  & \vdots  \\
   V_lV_1^*  & \cdots & |V_l|^2 & \cdots &  V_lV_n^*   \\
      \vdots  &  & \vdots & \ddots & \vdots  \\
   V_nV_1^* & \cdots & V_nV_l^* & \cdots &  |V_n|^2   \\
   \end{matrix}  \right]$$
with $X_{lm}$ the corresponding element in the $l$th row and $m$th column. By ignoring the rank constraints, the standard SDP relaxation of the OPF problem is obtained with the following set of constraints:
\begin{align}
\tag{SDP-OPF} \label{eq:sdpopf} \\
%\min\ \ &\sum_{l\in G} c^2_l(P^g_l)^2+c^1_l(P^g_l)+c^0_l \label{eq1}\\
&P_l^{\min} \leq P_l^g \leq P_l^{\max}, \ \forall l \in G,  \label{eq_2}\\
& Q_l^{\min}\leq Q^g_l\leq Q_l^{\max}, \ \forall l \in G, \label{eq_3}\\
& P^g_l -P_l^d= \sum_{m\in N(l)} P_{lm}, \ \forall l \in N\\
& Q^g_l - Q_l^d= \sum_{m\in N(l)} Q_{lm}, \ \forall l \in N  \\   
& (V_l^{\min})^2 \leq X_{ii} \leq  (V_l^{\max})^2,  \ \forall l \in N, \label{eq_4}\\
& P_{lm}^2+Q_{lm}^2 \leq (S_{lm}^{\max})^2,  \ \forall (l,m)\in E,\label{eq_5} \\
& P_{lm} =g_{lm} X_{ll}-g_{lm}\Re{X_{lm}}- b_{lm} \Im{X_{lm}} \\
&Q_{lm}  = -(b_{lm}+\frac{\bar{b}_{lm}}{2})  X_{ll} +b_{lm}\Re{X_{lm}} - g_{lm} \Im{X_{lm}}\\
& X=X^* \succeq 0
\end{align}
The SDP relaxation in the complex domain is formulated in Bose et al. \cite{Bose_2014} and is widely used in the literature now for its notational simplicity. The SDP relaxation for OPF has been widely used since it was originaly proposed by Bai et al. \cite{Bai_2008} and  Lavaei  and  Low  \cite{Lavaei_2012}. Because convex conic programs are polynomially solvable, the SDP relaxation offers an effective way for obtaining global optimal solutions to OPF problems whenever the relaxation is exact. Motivated by the success of SDP relaxations for OPF problems, an SDP approach is proposed to solve the TNEP problem.

TNEP consists of finding the minimum cost plan for the electrical system expansion so that the network adequately serves the forecasted system load over a given horizon; the OPF problem appears as the backbone of TNEP. In addition to the variables defined above, let
$$\alpha^t_{lm}=\begin{cases}
                1 & \text{ if line $t$ is installed in transmission corridor $(l,m)$ }\\
                0 & \text{ otherwise.}
               \end{cases}
$$
Two additional complex matrices are defined, $X^t_{lm}$ and $X_{l(lm)}^t$, the values of $X_{lm}$ and $X_{ll}$ for the $t$ line in transmission corridor $(l,m)$ respectively. When line $t$ is not installed in the transmission corridor $(l,m)$, the corresponding $X_{l(lm)}^t=0$ and the real and imaginary parts of $X^t_{lm}=0$. While when line $t$ is installed, the corresponding $X_{l(lm)}^t=X_{ll}=| V_l |^2$ and the real and imaginary parts of $X^t_{lm}$ are equal to that of $X_{lm}$ and are bounded by:
      $$0 \leq \Re\{X_{lm}^t\}\leq V_l^{\max}V_m^{\max}$$  
      $$-V_l^{\max}V_m^{\max} \leq \Im\{X_{lm}^t\}\leq V_l^{\max}V_m^{\max}.$$
Thus, the TNEP optimization model can be formulated as a mixed integer semidefinite program (MISDP) that is a relaxation of the mixed-integer ACOPF problem:
{\footnotesize
\begin{subequations}
   \begin{align}
   \tag{TNEP} \label{eq:tnep} \\
{\min \: \:} & f(\alpha,P^g_l)={\sum_{(l,m)\in E} \sum_{t=n^{\min}}^{n^{\max}}c_{lm}\alpha^t_{lm}+c_p \sum_{l=1}^n P^g_l} \label{eqn1} \\
    \text{s.t. }& P^g_l - P^d_l = \sum_{m\in N(l)} P_{lm} \label{eqn2} \\
		& Q^g_l - Q^d_l = \sum_{m\in N(l)} Q_{lm}\label{eqn3}  \\
                & P_l^{\min} \le P^g_l \le P_l^{\max}   \label{eqn4}\\
                & Q_l^{\min} \le Q^g_l \le Q_l^{\max} \label{eqn5}\\
                & P_{lm}=  \sum_{t=1}^{n_{lm}^{\max}}P^t_{lm}  \label{eqn6}\\
                & Q_{lm}=  \sum_{t=1}^{n_{lm}^{\max}}Q^t_{lm}  \label{eqn7}\\
                & P_{lm}^t=  g_{lm}X_{l(lm)}^t-g_{lm}\Re\{X_{lm}^t\}-b_{lm}\Im\{X_{lm}^t\}  \label{eqn8}\\
                & Q_{lm}^t=  -(b_{lm}+\frac{\bar{b}_{lm}}{2})X_{l(lm)}^t+b_{lm}\Re\{X_{lm}^t\}-g_{lm}\Im\{X_{lm}^t\}  \label{eqn9} \\
                & n_{lm}^{\min} \leq \sum_{t=1}^{n_{lm}^{\max}}\alpha^t_{lm}  \leq n_{lm}^{\max} \label{eqn10}\\
                 & \alpha^t_{lm} \leq \alpha^{t-1}_{lm}  \label{eqn10a}\\
                & (V_l^{\min})^2 \alpha^t_{lm} \leq X^t_{l(lm)} \leq (V_l^{\max})^2 \alpha^t_{lm} \label{eqn11}\\ 
                & 0 \leq \Re\{X_{lm}^t\}\leq V_l^{\max}V_m^{\max} \alpha^t_{lm}  \label{eqn12}\\ 
                 &-V_l^{\max}V_m^{\max} \alpha^t_{lm}   \leq \Im\{X_{lm}^t\}\leq V_l^{\max}V_m^{\max} \alpha^t_{lm} \label{eqn13}\\ 
                 &(V_l^{\min})^2(1- \alpha^t_{lm})   \leq X_{ll}- X_{l(lm)}^t \leq (V_l^{\max})^2(1- \alpha^t_{lm}) \label{eqn14} \\ 
                 &  0 \leq \Re\{X_{lm}\}- \Re\{X_{lm}^t\}\leq V_l^{\max}V_m^{\max} (1-\alpha^t_{lm}) \label{eqn15}\\ 
                  &   -V_l^{\max}V_m^{\max} (1-\alpha^t_{lm})  \leq \Im\{X_{lm}\}- \Im\{X_{lm}^t\} \label{eqn16}\\
                  &  \Im\{X_{lm}\}- \Im\{X_{lm}^t\} \leq V_l^{\max}V_m^{\max} (1-\alpha^t_{lm}) \label{eqn17}\\ 
                  & {(P_{lm}^1)^2+(Q_{lm}^1)^2} \leq (S_{lm}^{\max})^2\alpha_{lm}^1  \label{eqn18}\\
                &\alpha_{lm}^t \in \{0,1\},  \qquad   X =X^*\succeq 0.\label{eqn19}
\end{align}
\end{subequations}
}
The objective function \eqref{eqn1} minimizes the expansion cost with a penalty term on the power losses \cite{Jabr_2013}. Constraints \eqref{eqn2}-\eqref{eqn5} set limits on the active power and reactive power on each bus. Constraints~\eqref{eqn6}-\eqref{eqn9} compute the total real/reactive power flow along a transmission corridor which is obtained by summing the power flow over the distinct lines. Constraints~\eqref{eqn10} and \eqref{eqn10a} enforce a limit on the number of lines added and ensure a sequential installation of circuits in each transmission line respectively. Constraints~\eqref{eqn11}-\eqref{eqn17} enforce that when line $t$ is not installed in the transmission corridor $(l,m)$ then $\alpha^t_{lm}=0$ and the corresponding $X_{l(lm)}^t$=0. In addition, the real and imaginary parts of $X^t_{lm}=0$ and hence $P_{lm}^t= Q_{lm}^t =0$. While when line $t$ is installed then $\alpha^t_{lm}=1$ and the corresponding $X_{l(lm)}^t=X_{ll}=\mid V_l \mid^2$. In addition, the real and imaginary parts of $X^t_{lm}$ are equal to that of $X_{lm}$ and are bounded by the voltage constraints. Finally, constraints~\eqref{eqn18} set a limit on the apparent power flow at every line. The presence of the binary variables in addition to the SDP constraint in the TNEP formulation, \eqref{eq:tnep}, results in a mixed integer semidefinite program that is difficult to solve.

\section{Valid Inequalities}
To solve the above formulation of the TNEP problem, \eqref{eq:tnep}, a specialized branch-and-bound algorithm is developed where the binary condition is relaxed and at every node a SDP relaxation is solved. Since SDP relaxations are typically expensive to solve, having a strong bound at the root node is critical to the computational efficiency of the algorithm. Thus, a set of valid inequalities is introduced to attain a stronger relaxation of the feasible region of the TNEP problem. In  order  to  obtain  the  inequalities, results from \cite{Haffner2000} are used to generate fencing constraints. Additionally, conic constraints are added to improve the quality of the relaxation. The  resulting  convex  program is a strengthened  SDP program with improved lower bounds. %Hence, to strengthen the lower bound further valid inequalities are added to the SDP relaxation in a branch-and-cut framework.

The first set of constraints that are added are conic constraints:
\begin{subequations}
  \begin{align}
& {(P_{lm}^t)^2+(Q_{lm}^t)^2} \leq (S_{lm}^{\max})^2\alpha_{lm}^t  & \:& \forall \, t, \, \forall \,  (l,m) \in E.
 \end{align}
\end{subequations}   

The second set of constraints are referred to as the fence constraints. Fencing constraints are a generalization of Kirchhoff Current Law, KCL and form part of a heuristic methodology for transmission expansion planning \cite{Haffner2000}. Three kinds of fences are generated and added to the SDP relaxation: around a single node, around one node and neighboring node, and around a node and its entire neighborhood. To generate a fencing constraint:
\begin{itemize}
 \item Place an imaginary fence around a portion of the power system.
\item Calculate power delivered into or out of the fenced area.
\item Compare the transmission capacity with the net load or net generation within the fenced area to determine if the transmission is adequate.
\end{itemize}
The constraints are written as:
\begin{subequations}
  \begin{align}
&         \bar{n}_m=\min\{n^{\max}, \ceil[\bigg]{\frac{\sum_lP_l^d-\sum_lP_l^{\max}}{S_{lm}^{\max}}}\} \\
&         \sum_{m\in N(l)}\sum_{t=n^{\min}}^{\bar{n}_m}\alpha^t_{lm} \geq \ceil[\bigg]{\frac{\sum_l P_l^d-\sum_l P_l^{\max}}{S^{\max}}}
 \end{align}
\end{subequations} 
where $\lceil . \rceil$ is the ceiling operator and $S^{\max}=\max_{(l,m)} \{S_{lm}^{\max}\}$. 
\begin{example}
Consider the example with 6 nodes as shown in Figure \ref{fig1} \cite{Rider_2007}. A fence is set on nodes 1-5 of the power network. The sum of the loads at buses 1-5 is $\sum_l P^d_l$ =760 while the sum of max power generated at the same buses is $\sum_l P_l^{\max}$=530 and $n^{\max}=5$. The power flow limitation is given as $S_{16}^{\max}=90$, $S_{56}^{\max}=98$, $S_{l6}^{\max}=120$ $(l=2,3,4)$. The transmission  deficit can be represented by the following fencing constraint \ \\ $\alpha^1_{1,6}+\alpha^1_{2,6}+\alpha^1_{3,6}+\alpha^1_{4,6}+\alpha^1_{5,6}+\alpha^2_{1,6}+\alpha^2_{2,6}+\alpha^2_{3,6}+\alpha^2_{4,6}+\alpha^2_{5,6}+\alpha^3_{1,6}+\alpha^3_{5,6} \geq 2.$
\begin{figure}
\centering
{\includegraphics[trim=0cm 0cm 0cm 0cm, clip=true, scale=0.45]{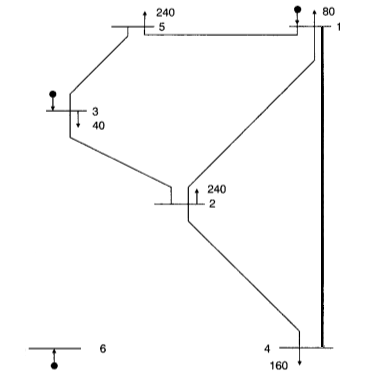}}
\caption{Example of a 6 bus network.} \label{fig1}
\end{figure}

\end{example}

\section{Sparsity of the SDP Relaxation}
Exploiting sparsity has been one of the essential tools for solving large-scale optimization problems in general and optimal power flow problems in particular \cite{Ghaddar2016, Molzahn2015}. Solving existing semidefinite relaxations of large optimal power flow problems requires exploiting power system sparsity. Using a matrix completion decomposition, existing semidefinite relaxations of optimal power flow problems are computationally tractable for problems with hundreds of buses. Since the optimal power flow problem is a building block of transmission expansion network planning, to solve the latter a framework for exploiting the sparsity characterized in terms of a chordal graph structure via positive semidefinite matrix completion is used \cite{Grone1984}. 

The computation time of solving the SDP relaxation can be significantly improved using SparseColO that is used as a preprocessor, which reduces the dimension of matrix variables in an SDP relaxation before applying the SDP solver. When applied to the semidefinite relaxation, SparseColO enhances the structured sparsity of the problem, i.e., the correlative sparsity. As a result, the resulting SDP can be solved more effectively by applying the sparse SDP relaxation. In a branch-and-cut framework, exploiting sparsity is important as it reduces the computational time at each node of the tree resulting in a significant decrease in the total solving time.
\section{SDP Branch-and-Cut}
To develop an efficient technique, we adopt a branch-and-bound approach to solve the TNEP problem using a SDP relaxation, combined with valid inequalities. The success and the computational efficiency of the branch-and-bound procedure strongly depends on the quality of the relaxation bounds, the early generation of binary feasible solutions to get good bounds, and the branching rules used to obtain the subproblems \cite{NW1988}. Applying the SDP relaxation in addition to the valid inequalities can help speed up the branch-and-bound process by improving the bounds at each node, thus reducing the number of nodes of the tree. The following sections, describe the branch-and-bound algorithm that is implemented to solve TNEP. The general algorithm, feasible solution, and branching rules are described in details. 
\subsection{Feasible Solution}\label{sec-feasible}
During a branch-and-bound procedure, it is of particular importance to find a feasible solution to be able to fathom nodes in the tree reducing the search space. In the described branch-and-bound algorithm, rounding the optimal solution $\alpha^t_{lm}$ to obtain $\hat{\alpha}^t_{lm}$ is utilized for a feasible solution. {The obtained rounded solution is checked for feasibility and the rank of the SDP matrix is checked. The incumbent is updated in case the solution is feasible and the rank of the SDP matrix is 1.} For the branching strategy, depth first is used, that is, the next node to be solved of the branch-and-bound tree is one of the child nodes of the current node solved. Depth-first node selection goes deep into the branch and bound tree at each iteration, so it reaches the leaf nodes quickly. This is one way of achieving an early feasible solution and hence an incumbent. 
\subsection{Inequality Generation Scheme}
Branch-and-cut methods combine the classical branch-and-bound with cutting plane techniques; they define valid linear inequalities to improve the lower bound of the solution and speed up the search. The branch-and-cut approach can be sped up considerably by the employment of an inequality generation scheme, either just at the root node of the tree, or at every node of the tree. In the case of the TNEP problem, the inequalities are only generated at the root node, but since they are valid to all the children nodes, they are added to each node of the tree. 
\subsection{Branch-and-Cut Algorithm}\label{sec-algorithm}
As discussed in the previous sections, the branch-and-cut algorithm for TNEP problem consists of solving SDP relaxations at every node of the tree. Since the variables are binary, the set of possible configurations of the variables is finite and equal to $2^{n^{\max}*|E|}$ thus the branch-and-bound procedure stops after a finite number of nodes. Given (TNEP) and the binary index set $J\subseteq \{1,\ldots,n^{\max}\}$, a sketch of the branch-and-cut algorithm is given as follows:
\begin{algorithm}[ht]\caption{Branch-and-Cut Algorithm for TNEP}\label{alg:bandb}
\textbf{input}: {(TNEP)  }\BlankLine
\textbf{output}: {$\alpha^t_{lm}$ and $z_{TNEP}$}\BlankLine
\textbf{set}:{$\:G_0=G\cup \{ \alpha^t_{lm} \in \{0,1\}: t \in J \}, \: \mu_{G_0} = -\infty$, $z_{TNEP}= \infty$, $\hat{\alpha}^t_{lm}=[\:]$, $Nodes=\{N_{0}\}$}
\BlankLine
 \While{$Nodes \neq \phi$}{
 Choose node $N_i \in Nodes$\;
 \textbf{set}: $Nodes=Nodes \backslash N_i$ \;
 Solve (TNEP$_i$) to obtain $z_i$ and ${\alpha}^t_{lm}$\;
 \uIf{$z_i=-\infty$ }{ Fathom $N_i$\;}
 \uElseIf {$z_i \geq  z_{TNEP}$,} {Fathom $N_i$\;}
  \uElseIf {${\alpha}^t_{lm} \in \{0,1\} \: \forall t \: \forall (l,m)$,} { 
  \If{$f({\alpha}^t_{lm}) \leq z_{TNEP}$}{$z_{TNEP} = f({\alpha}^t_{lm}), \: \hat{\alpha}^t_{lm}={\alpha}^t_{lm}$\;}
  Fathom $N_i$\;}
\Else {Apply rounding to ${\alpha}^t_{lm}$ to obtain $\bar{\alpha}^t_{lm}$\;
 \If{feasible and $f(\bar{\alpha}^t_{lm}) \leq z_{TNEP}$}{$z_{TNEP} = f(\bar{\alpha}^t_{lm}), \: \hat{\alpha}^t_{lm}=\bar{\alpha}^t_{lm}$\;}
 Choose $(G_1,G_2)$ as a branching rule for $G$\;}
 \textbf{set}: $ Nodes = Nodes \cup \{ N_{G_1}, N_{G_2}\}$\;}
\end{algorithm}
%\eIf{$\mu_G =f(\bar{\alpha}^t_{lm})$}{
% Fathom $N_i$\;}
Valid inequalities are added at the root node and are reused at the nodes of the branch-and-bound tree as they are also valid for the children nodes.
\section{Computational Results}
In this section, the performance of the proposed branch-and-cut approach is evaluated for the TNEP problem on different test cases. The method was programmed in MATLAB and the computational experiments were conducted on a Lenovo Thinkstation P300 with 32GB of RAM. The MOSEK \cite{mosek} solver was used to solve the SDP relaxations and SparseColO \cite{kim2011} was used for exploiting sparsity  at each node of the branch-and-bound tree. Four networks are considered, two 6-bus \cite{Rider_2007}, 24-bus \cite{Rider_2007}, and 46-bus systems \cite{Jabr_2013}. The 6-bus Graver system has six buses and 75 potential lines for the addition of new circuits. From the basic topology, two scenarios are taken into account: expansion with existing lines (instance 6-bus), and expansion without existing lines (instance 6gf-bus). The third instance is the IEEE 24-bus and the fourth test system is the 46-bus south Brazilian network.

Table \ref{tb:results}, presents the computational results for the proposed branch-and-bound approach for the SDP relaxation with and without valid inequalities as well as with and without exploiting sparsity. The root node lower bound as well as the optimal solution (scaled objective value \eqref{eqn1}) are reported in the table. Additionally, the computational time (in seconds) and the number of nodes of the branch-and-bound tree are given.
From Table \ref{tb:results}, it can be observed that the strengthened SDP relaxation with valid inequalities has a stronger lower bound at the root node and this translates to lower number of nodes as well as computational time. At the end of the branch-and-bound algorithm, the optimal solution is obtained for the original non-convex quadratic problems as the SDP relaxation produces a rank-1 solution. Notice the significant improvement in the root bound for 6gf instance where the root bound improved from 55.02 to 145.98. Furthermore, the computational time of solving TNEP using sparsity is significantly lower than the computational time of solving the problem without exploiting sparsity. The solution time of the SDP relaxation of the instances with 24 and 46 buses are reduced by a factor of around 10 when exploiting sparsity.
\begin{table*}[!htb]
 \centering
\caption{Computational results for the branch-and-bound approach with and without valid inequalities and with and without exploiting sparsity.}
\begin{tabular}{l|rrrr||rrrr||rrrr}
 &\multicolumn{4}{c||}{MISDP B\&B}&\multicolumn{4}{c}{MISDP + Valid Inequalities B\&C}&\multicolumn{4}{c}{MISDP + Valid Inequalities + Sparsity}\\
$|N|$&Root & Optimal& Nodes & CPU & Root & Optimal& Nodes&CPU& Root & Optimal& Nodes&CPU\\
\hline\hline
6   & 29.04 &   160.12 & 1077 & 527.71   & 77.80  &160.12 & 665 & 358.12& 77.80  &160.12 & 665 & 224.09 \\
6gf & 55.02 &	250.11 & 6669 & 2707.01 & 145.98 & 250.11 & 2515& 1101.30& 145.98 & 250.11 & 2515& 891.48 \\
24  & 7.38  &	88.29  &1133  & 34430.00 & 44.17 & 88.29  & 337 & 10028.51 & 44.17 & 88.29  & 337 & 955.81 \\
46  & 7.88  & 	-      & -    & -        & 37.69  & 74.75  & 1603 & 90108.03 & 37.69  & 74.75  & 1603 & 8191.64  \\
\hline
\end{tabular}\label{tb:results}
\end{table*}
Table \ref{tb:data} summarizes the number of transmission corridors (TC) and the total number of existing (EL), potential (PL), and optimal lines (OL) for all test systems. The expansion plans for the four test systems are given in Tables \ref{tb:results1}-\ref{tb:results3}, together with their expansion costs. By comparing with the results in \cite{Jabr_2013}, it becomes evident that the mixed-integer conic solver \cite{Jabr_2013} did give a globally optimal solution for the same TNEP instances although it could not rigorously prove the optimality of the designs.

Also as in \cite{Jabr_2013}, a comparison with the scaled objective values in Table \ref{tb:results} shows that the losses are a small fraction of the total objective \eqref{eqn1} and are therefore not expected to affect the optimal design; in fact, the loss term is included in the objective function so as to induce a rank-1 solution.

\begin{table}[!h]
\centering
\caption{Test System and Optimal Data}
\begin{tabular}{l|rrrrr}
$|N|$&TC & EL & PL & OL\\
\hline\hline
6   & 15& 6& 75& 6 \\
6gf & 15& 0& 75& 10 \\
24  & 41& 38& 123& 3 \\
46  & 79& 62& 158& 8\\
\hline
\end{tabular}\label{tb:data}
\end{table}

\begin{table}[!h]
 \centering
\caption{Number of circuits and total cost produced by MISDP for the 6-bus network.}
\begin{tabular}{l|r||r}
 &\multicolumn{1}{l||}{initial }&\multicolumn{1}{l}{without initial }\\
$i$-$j$  &\multicolumn{1}{l||}{ network}&\multicolumn{1}{l}{network}\\
\hline\hline
1-5  &   0 & 1\\
2-3 &   0 & 1 \\
2-6  &  2 & 2 \\
3-5  & 2  &3  \\
4-6&   2  & 3 \\ \hline
cost ($\times10^3$\$) &160 & 250 \\
\hline
\end{tabular}\label{tb:results1}
\end{table}

\begin{table}[!h]
 \centering
\caption{Number of circuits and total cost produced by MISDP for the 24-bus network.}
\begin{tabular}{l|r}
$i$-$j$ & MISDP \\
\hline\hline
6-10  &  1 \\
7-8 & 1  \\
14-16&  1  \\ \hline
cost ($\times10^3$\$) & 86   \\
\hline
\end{tabular}\label{tb:results2}
\end{table}
\begin{table}[!h]
 \centering
\caption{Number of circuits and total cost produced by MISDP for the 46-bus network.}
\begin{tabular}{l|rr}
$i$-$j$ & MISDP \\
\hline\hline
5-6  &  2 \\
19-25 &  0  \\
20-21   & 2   \\
20-23    & 1   \\
24-25&   0  \\ 
33-34&   1  \\ 
42-43&   1  \\ 
42-44&   0  \\ 
46-6&  1  \\ \hline
cost ($\times10^3$\$) &71,451   \\
\hline
\end{tabular}\label{tb:results3}
\end{table}

\section{Conclusion}
This paper presented the implementation details of a mixed-integer
SDP method that is specialized to solve the AC-TNEP problem. The method
is demonstrated on standard TNEP problems with 6, 24, and 46 nodes.
Valid inequalities and exploiting the problem structure through sparsity are proposed to speed up the branch-and-bound search;
without these cuts, the solution of the largest network is not even
possible. The planning results are feasible with respect to the AC
power flow constraints, and they are globally optimal; this is in
contrast to many of the results that are obtained from simplified
linear models, and which are not even feasible from an AC network
standpoint. Although previous AC-TNEP methods reported some test results
that are identical to the ones herein, the importance of this work
is that it is the first to rigorously ascertain global optimality.
This has value in checking the quality of heuristic and local search
approaches, in addition to solutions that are based on simplified
modeling.

Future directions include extending the model to deal with potential changes in the demand as well as uncertainties in the renewable power generation in transmission expansion planning. Furthermore, the same approach can be applied to operational problems such as operations with transmission switching.

\section*{Acknowledgment}
Bissan Ghaddar was supported by NSERC Discovery Grant RGPIN-2017-04185.
\bibliographystyle{IEEEtran}
\bibliography{TNEP_refs}

\begin{thebibliography}{10}
\providecommand{\url}[1]{#1}
\csname url@samestyle\endcsname
\providecommand{\newblock}{\relax}
\providecommand{\bibinfo}[2]{#2}
\providecommand{\BIBentrySTDinterwordspacing}{\spaceskip=0pt\relax}
\providecommand{\BIBentryALTinterwordstretchfactor}{4}
\providecommand{\BIBentryALTinterwordspacing}{\spaceskip=\fontdimen2\font plus
\BIBentryALTinterwordstretchfactor\fontdimen3\font minus
  \fontdimen4\font\relax}
\providecommand{\BIBforeignlanguage}[2]{{%
\expandafter\ifx\csname l@#1\endcsname\relax
\typeout{** WARNING: IEEEtran.bst: No hyphenation pattern has been}%
\typeout{** loaded for the language `#1'. Using the pattern for}%
\typeout{** the default language instead.}%
\else
\language=\csname l@#1\endcsname
\fi
#2}}
\providecommand{\BIBdecl}{\relax}
\BIBdecl

\bibitem{Garver_1970}
L.~L. Garver, ``Transmission network estimation using linear programming,''
  \emph{IEEE Trans. Power Appar. Syst.}, vol. PAS-89, no.~7, pp. 1688--1697,
  Sep. 1970.

\bibitem{Escobar_2004}
A.~H. Escobar, R.~A. Gallego, and R.~Romero, ``Multistage and coordinated
  planning of the expansion of transmission systems,'' \emph{IEEE Trans. Power
  Syst.}, vol.~19, no.~2, pp. 735--744, May 2004.

\bibitem{Vinasco_2011}
G.~Vinasco, M.~J. Rider, and R.~Romero, ``A strategy to solve the multistage
  transmission expansion planning problem,'' \emph{IEEE Trans. Power Syst.},
  vol.~26, no.~4, pp. 2574--2576, Nov. 2011.

\bibitem{Romero_2002}
R.~Romero, A.~Monticelli, A.~Garcia, and S.~Haffner, ``Test systems and
  mathematical models for transmission network expansion planning,'' \emph{IEE
  Proc. - Gener. Transm. Distrib.}, vol. 149, no.~1, pp. 27--36, Jan. 2002.

\bibitem{Lavaei_2012}
J.~Lavaei and S.~H. Low, ``Zero duality gap in optimal power flow problem,''
  \emph{IEEE Trans. Power Syst.}, vol.~27, no.~1, pp. 92--107, Feb. 2012.

\bibitem{Haffner_2000}
S.~Haffner, A.~Monticelli, A.~Garcia, J.~Mantovani, and R.~Romero, ``Branch and
  bound algorithm for transmission system expansion planning using a
  transportation model,'' \emph{IEE Proc. - Gener. Transm. Distrib.}, vol. 147,
  no.~3, pp. 149--156, May 2000.

\bibitem{Haffner_2001}
S.~Haffner, A.~Monticelli, A.~Garcia, and R.~Romero, ``Specialised
  branch-and-bound algorithm for transmission network expansion planning,''
  \emph{IEE Proc. - Gener. Transm. Distrib.}, vol. 148, no.~5, pp. 482--488,
  Sep. 2001.

\bibitem{Villasana_1985}
R.~Villasana, L.~L. Garver, and S.~J. Salon, ``Transmission network planning
  using linear programming,'' \emph{IEEE Trans. Power Appar. Syst.}, vol.
  PAS-104, no.~2, pp. 349--356, Feb. 1985.

\bibitem{Romero_2007}
R.~Romero, E.~N. Asada, E.~Carreno, and C.~Rocha, ``Constructive heuristic
  algorithm in branch-and-bound structure applied to transmission network
  expansion planning,'' \emph{IET Gener. Transm. Distrib.}, vol.~1, no.~2, pp.
  318--323, Mar. 2007.

\bibitem{Rider_2007}
M.~J. Rider, A.~V. Garcia, and R.~Romero, ``Power system transmission network
  expansion planning using {AC} model,'' \emph{IET Gener. Transm. Distrib.},
  vol.~1, no.~5, pp. 731--742, Sep. 2007.

\bibitem{Bahiense_2000}
L.~Bahiense, G.~C. Oliveira, M.~Pereira, and S.~Granville, ``A mixed integer
  disjunctive model for transmission network expansion,'' \emph{IEEE Trans.
  Power Syst.}, vol.~16, no.~3, pp. 560--565, Aug. 2001.

\bibitem{Zhang_2012_c}
H.~Zhang, G.~T. Heydt, V.~Vittal, and H.~D. Mittelmann, ``Transmission
  expansion planning using an {AC} model: Formulations and possible
  relaxations,'' in \emph{IEEE Power and Energy Society General Meeting}, Jul.
  2012, pp. 1--8.

\bibitem{Bent_2012}
R.~Bent, G.~L. Toole, and A.~Berscheid, ``Transmission network expansion
  planning with complex power flow models,'' \emph{IEEE Trans. Power Syst.},
  vol.~27, no.~2, pp. 904--912, May 2012.

\bibitem{Gallego_1998a}
R.~A. Gallego, A.~Monticelli, and R.~Romero, ``Transmission system expansion
  planning by an extended genetic algorithm,'' \emph{IEE Proc. - Gener. Transm.
  Distrib.}, vol. 145, no.~3, pp. 329--335, May 1998.

\bibitem{Gallego_1998b}
------, ``Comparative studies on nonconvex optimization methods for
  transmission network expansion planning,'' \emph{IEEE Trans. Power Syst.},
  vol.~13, no.~3, pp. 822--828, Aug. 1998.

\bibitem{da_Silva_2000}
E.~L. {da Silva}, H.~A. Gil, and J.~M. Areiza, ``Transmission network expansion
  planning under an improved genetic algorithm,'' \emph{IEEE Trans. Power
  Syst.}, vol.~15, no.~3, pp. 1168--1174, Aug. 2000.

\bibitem{da_Silva_2001}
E.~L. {da Silva}, J.~M.~A. Ortiz, G.~C. {de Oliveira}, and S.~Binato,
  ``Transmission network expansion planning under a tabu search approach,''
  \emph{IEEE Trans. Power Syst.}, vol.~16, no.~1, pp. 62--68, Feb. 2001.

\bibitem{de_J_Silva_2005}
I.~{de J Silva}, M.~J. Rider, R.~Romero, A.~V. Garcia, and C.~A. Murari,
  ``Transmission network expansion planning with security constraints,''
  \emph{IEE Proc. - Gener. Transm. Distrib.}, vol. 152, no.~6, pp. 828--836,
  Nov. 2005.

\bibitem{Rodriguez_2008}
J.~I.~R. Rodriguez, D.~M.~F. {a}o, and G.~N. Taranto, ``Short-term tranmission
  expansion planning with {AC} network model and security constraints,'' in
  \emph{Power Systems Computation Conference}, Jul. 2008, pp. 1--7.

\bibitem{Rahmani20101056}
M.~Rahmani, M.~Rashidinejad, E.~Carreno, and R.~Romero, ``Efficient method for
  {AC} transmission network expansion planning,'' \emph{Elec. Power Syst.
  Res.}, vol.~80, no.~9, pp. 1056 -- 1064, 2010.

\bibitem{Romero_1994a}
R.~Romero and A.~Monticelli, ``A hierarchical decomposition approach for
  transmission network expansion planning,'' \emph{IEEE Trans. Power Syst.},
  vol.~9, no.~1, pp. 373--380, Feb. 1994.

\bibitem{Romero_1994b}
------, ``A zero-one implicit enumeration method for optimizing investments in
  transmission expansion planning,'' \emph{IEEE Trans. Power Syst.}, vol.~9,
  no.~3, pp. 1385--1391, Aug. 1994.

\bibitem{Binato_2001}
S.~Binato, M.~V.~F. Pereira, and S.~Granville, ``A new benders decomposition
  approach to solve power transmission network design problems,'' \emph{IEEE
  Trans. Power Syst.}, vol.~16, no.~2, pp. 235--240, May 2001.

\bibitem{AlizadehMousavi2016275}
O.~Alizadeh-Mousavi and M.~Zima-Bo\v{c}karjova, ``Efficient benders cuts for
  transmission expansion planning,'' \emph{Elec. Power Syst. Res.}, vol. 131,
  pp. 275 -- 284, 2016.

\bibitem{Zhang_2012}
H.~Zhang, V.~Vittal, G.~T. Heydt, and J.~Quintero, ``A mixed-integer linear
  programming approach for multi-stage security-constrained transmission
  expansion planning,'' \emph{IEEE Trans. Power Syst.}, vol.~27, no.~2, pp.
  1125--1133, May 2012.

\bibitem{Zhang_2013}
H.~Zhang, G.~T. Heydt, V.~Vittal, and J.~Quintero, ``An improved network model
  for transmission expansion planning considering reactive power and network
  losses,'' \emph{IEEE Trans. Power Syst.}, vol.~28, no.~3, pp. 3471--3479,
  Aug. 2013.

\bibitem{Bent_2014}
R.~Bent, C.~Coffrin, R.~R. {Gumucio E.}, and P.~V. Hentenryck, ``Transmission
  network expansion planning: Bridging the gap between {AC} heuristics and {DC}
  approximations,'' in \emph{Power Systems Computation Conference}, Aug. 2014,
  pp. 1--8.

\bibitem{Macedo_2016}
L.~H. Macedo, C.~V. Montes, J.~F. Franco, M.~J. Rider, and R.~Romero, ``{MILP}
  branch flow model for concurrent {AC} multistage transmission expansion and
  reactive power planning with security constraints,'' \emph{IET Gener. Transm.
  Distrib.}, vol.~10, no.~12, pp. 3023--3032, 2016.

\bibitem{Dominguez_2017}
A.~H. Dominguez, L.~H. Macedo, A.~H. Escobar, and R.~Romero, ``Multistage
  security-constrained {HVAC/HVDC} transmission expansion planning with a
  reduced search space,'' \emph{IEEE Trans. Power Syst.}, vol.~PP, no.~99, pp.
  1--1, 2017.

\bibitem{Jabr_2013}
R.~A. Jabr, ``Optimization of {AC} transmission system planning,'' \emph{IEEE
  Trans. Power Syst.}, vol.~28, no.~3, pp. 2779--2787, Aug. 2013.

\bibitem{Bose_2014}
S.~Bose, S.~H. Low, T.~Teeraratkul, and B.~Hassibi, ``Equivalent relaxations of
  optimal power flow,'' CALTECH, Tech. Rep., 2014, available at
  \url{http://arxiv.org/pdf/1401.1876v1.pdf}.

\bibitem{Bai_2008}
X.~Bai, H.~Wei, K.~Fujisawa, and Y.~Wang, ``Semidefinite programming for
  optimal power flow problems,'' \emph{Int. J. Elec. Power}, vol.~30, no.~6,
  pp. 383--392, 2008.

\bibitem{Haffner2000}
S.~Haffner, A.~Monticelli, A.~Garcia, J.~Mantovani, and R.~Romero, ``Branch and
  bound algorithm for transmission system expansion planning using a
  transportation model,'' \emph{IEE Proceedings - Generation, Transmission and
  Distribution}, vol. 147, no.~3, pp. 149--156, 2000.

\bibitem{Ghaddar2016}
B.~Ghaddar, J.~Marecek, and M.~Mevissen, ``Optimal power flow as a polynomial
  optimization problem,'' \emph{IEEE Transactions on Power Systems}, vol.~31,
  no.~1, pp. 539--546, 2016.

\bibitem{Molzahn2015}
D.~K. Molzahn and I.~A. Hiskens, ``Sparsity-exploiting moment-based relaxations
  of the optimal power flow problem,'' \emph{IEEE Transactions on Power
  Systems}, vol.~30, no.~6, pp. 3168--3180, 2015.

\bibitem{Grone1984}
R.~Grone, C.~R. Johnson, E.~M. Sá, and H.~Wolkowicz, ``Positive definite
  completions of partial hermitian matrices,'' \emph{Linear Algebra and its
  Applications}, vol.~58, pp. 109 -- 124, 1984.

\bibitem{NW1988}
G.~Nemhauser and L.~Wolsey, \emph{Integer and Combinatorial
  Optimization}.\hskip 1em plus 0.5em minus 0.4em\relax Wiley, 1988.

\bibitem{mosek}
``The mosek optimization software,'' 2010, available at
  \url{http://www.mosek.com}.

\bibitem{kim2011}
S.~Kim, M.~Kojima, M.~Mevissen, and M.~Yamashita, ``Exploiting sparsity in
  linear and nonlinear matrix inequalities via positive semidefinite matrix
  completion,'' \emph{Mathematical Programming}, vol. 129, no.~1, pp. 33--68,
  2011.

\end{thebibliography}

\end{document}